\documentclass[leqno,11pt]{article}

\usepackage{amsmath}
\usepackage{epsfig}
\usepackage[latin1]{inputenc}
\usepackage{amssymb}
\usepackage[all]{xy}
\usepackage{color}

\newcommand{\qed}{\penalty 500\hfill$\square$\par\medskip}

\newtheorem{teor}{Theorem}

\newtheorem{lema}{Lemma}

\newtheorem{rem}{Remark}

\newcommand{\n}{\noindent}

\begin{document}
	
	\title{\Large{\bf Principal curves to fractional $m$-Laplacian systems and related
maximum and comparison principles}}
	
	\author{\sc
		A. L. A. de Araujo
		\thanks{Departamento de Matem\'{a}tica,
			Universidade Federal de Viçosa, 36570-900 - Vi\c{c}osa - MG, Brazil. E-mail: anderson.araujo@ufv.br}
		\,\,\,
		E. J. F. Leite
		\thanks{Departamento de Matem\'{a}tica,
			Universidade Federal de Viçosa, 36570-900 - Vi\c{c}osa - MG, Brazil. E-mail: edirjrleite@ufv.br}
		\,\,\,
		Aldo H. S. Medeiros
		\thanks{Departamento de Matem\'{a}tica,
			Universidade Federal de Viçosa, 36570-900 - Vi\c{c}osa - MG, Brazil. E-mail: aldo.medeiros@ufv.br}
	}
	
	\date{}

	\maketitle

	\setcounter{equation}{0}
	\begin{abstract}
In this paper we develop a comprehensive study on principal eigenvalues and both the (weak
and strong) maximum and comparison principles related to an important class of nonlinear systems involving fractional $m$-Laplacian operators.  Explicit lower bounds for principal eigenvalues of this system in terms of the diameter of $\Omega$ are also proved. As application, given $\lambda,\mu\geq 0$ we measure explicitly how small has to be $\text{diam}(\Omega)$ so that weak and strong maximum principles associated to this problem hold in $\Omega$.
		
	\end{abstract}
	
	\
	
	\
	
	\
	
	\date{}

	\maketitle

	\section{Introduction and main results}
\mbox{}	

The first part of the present work is dedicated to the study of principal eigenvalues for the following system of quasilinear elliptic boundary-value problems
	
	\begin{equation} \label{1}
		\left\{
		\begin{array}{llll}
			(-\Delta_p)^{s_1}u = \lambda a(x)\vert u \vert^{\alpha_1} \vert v \vert^{\beta_1 -1}v & {\rm in} \ \ \Omega,\\
			(-\Delta_q)^{s_2}v = \mu b(x)\vert v \vert^{\alpha_2} \vert u \vert^{\beta_2 -1}u & {\rm in} \ \ \Omega,\\
			u= v=0 & {\rm in} \ \mathbb{R}^{N}\setminus\Omega,
		\end{array}
		\right.
	\end{equation}
	where $\Omega$ is a bounded domain in $\mathbb{R}^N$ whose boundary is a with $C^{2}$-manifold, $p, q \in (1,\infty)$ and $s_1,s_2 \in (0,1)$ are given numbers, $a,b \in L^{\infty}(\Omega)$ are given functions satisfying
	\[
\underset{x\in\Omega}{\mathrm{ess }\inf}\ a(x)>0\ \ \text{ and }\ \ \underset{x\in\Omega}{\mathrm{ess }\inf}\ b(x)>0
\]
	and $\alpha_i, \beta_i$ are constants with $\alpha_i \geq 0$ and $\beta_i >0$ for $i = 1,2$.
	
	The fractional $m$-Laplacian operator $(-\Delta_m)^{s}$ is defined as
	
	\[
	(-\Delta_m)^{s}u(x) = C(N,s,m)\lim_{\varepsilon\searrow 0}\int\limits_{\mathbb{R}^{N}\setminus B_\varepsilon(x)}\frac{\vert u(x) - u(y) \vert^{m-2} (u(x)-u(y))}{\vert x-y\vert^{N+sm}}\; dy\, ,
	\]
	\n for all $x \in \mathbb{R}^{N}$, where $C(N,s,m)$ is a normalization factor. The fractional $m$-Laplacian is a nonlocal version of the $m$-Laplacian and is an extension of the fractional Laplacian $(m = 2)$. 
	
	These kinds of nonlocal operators have their applications in real world problems, such as obstacle problems, the study of American options in finance, game theory, image processing, and anomalous diffusion phenomena (see \cite{DiNezza} for more details). Due to this reason, elliptic problems involving the fractional Laplacian have been extensively studied in the two last decades.
	
	Throughout this paper to simplify notation, we omit the constant $C(N,s,m)$. The fractional Sobolev spaces $W^{s,m}(\Omega)$ is defined to be the set of functions $u \in L^m(\Omega)$ such that
	\[
	[u]_{s,m} = \left(\int\limits_{\mathbb{R}^{N}}\int\limits_{\mathbb{R}^{N}}\frac{\vert u(x) - u(y) \vert^{m}}{\vert x-y\vert^{N+sm}}\; dxdy\right)^{\frac{1}{m}} < \infty
	\]
	and we defined the space $W^{s,m}_0(\Omega)$ by
	\[
	W^{s,m}_0(\Omega) = \bigg\{ u \in W^{s,m}(\Omega); \ \ u = 0 \ \ \text{in} \ \ \Omega^c= \mathbb{R}^N \backslash \Omega\bigg\}
	\]
	this being uniformly convex, separable Banach space with norm $\Vert u \Vert := [u]_{s,m}$, whose
	dual space is denoted by $W^{-s,m'}(\Omega)$, where $m'=m/(m-1)$ (see \cite{DiNezza}). The embedding $W^{s,m}(\Omega) \hookrightarrow L^{\nu}(\Omega)$ is continuous for all $\nu \in [1,m^*_{s}]$ and compact for all $\nu \in [1,m^{*}_{s})$, with $m^{*}_{s} = \displaystyle\frac{Nm}{N-sm}.$
		
	Given any $f \in L^{\infty}(\Omega)$, we denote by $T_m^{s}(f): = u \in W^{s,m}_0(\Omega)$ the unique weak solution of the boundary-value problem
\begin{equation}\label{p1}
\left\{
\begin{array}{rrll}
(-\Delta_m)^s u &=& f(x) & {\rm in} \ \ \Omega,\\
u &=& 0  & {\rm in} \ \ \mathbb{R}^N\setminus \Omega.
\end{array}
\right.
\end{equation}
	
	It is well known that $u \in C^{\alpha}(\overline{\Omega})$ for some $0 < \alpha \leq s$ (see \cite{Iannizzotto}). Consider, for any $x \in \mathbb{R}^N$, the distance function, $d_{\Omega}: =  \text{dist}(x, \mathbb{R}^N \backslash \Omega)$ and the space
	\[
	C^0_s(\overline{\Omega}) =\bigg\{ u \in C^0(\overline{\Omega}); \frac{u}{d_{\Omega}^s} \ \text{has a continuos extension to }  \overline{\Omega}\bigg\}.
	\]
	
	Let also, $X := C^0_{s_1}(\overline{\Omega})\times C^0_{s_2}(\overline{\Omega})$, $X_+ := \bigg\{(f,g) \in X; \ \ f \geq 0 \ \ \text{and} \ \ g \geq 0 \ \ \text{in} \ \ \Omega \bigg\}$ and $\text{int}(X_+)$ is the topological interior of $X_+$ in $X$.

Let $u\in W^{s,m}_0(\Omega)$ be a weak solution of the problem (\ref{p1}). Then, by Corollary 4.2 of \cite{IMS}, we have
\begin{equation}\label{abp}
|| u||_{L^\infty(\Omega)}\leq C_{s,m}d^{\frac{ms}{m-1}}|| f||^{\frac{1}{m-1}}_{L^\infty(\Omega)},
\end{equation}
where $d:=\text{diam}(\Omega)$, $$C_{s,m}:=\sup\{w(x);x\in B_1\}$$
and $w\in W^{s,m}_0(B_1)\cap L^\infty(\mathbb{R}^N)$ is the unique weak solution of the problem

\[
\left\{
\begin{array}{rrll}
(-\Delta_m)^s w &= & 1 & {\rm in} \ \ B_1,\\
w &=&0  & {\rm in} \ \ \mathbb{R}^N\setminus B_1,
\end{array}
\right.
\]
where $B_1$ is the unit ball of $\mathbb{R}^N$.
	
	We say that $(\lambda, \mu) \in (\mathbb{R}_+^*)^2$ is an eigenvalue of the system \eqref{1} if possesses a nontrivial weak solution $(\varphi,\psi)$ in $W^{s_1,p}_0(\Omega)\times W^{s_2,q}_0(\Omega)$ which is called an eigenfunction corresponding to $(\lambda, \mu)$. We also say that $(\lambda, \mu)$ is a principal eigenvalue if possesses a positive eigenfunction $(\varphi,\psi)$; i.e., $\varphi$ and $\psi$ are positive in $\Omega$. Now, we say that $(\lambda, \mu)$ is simple in $\text{int}(X_+)$ if for any eigenfunctions $(\varphi,\psi),(\tilde{\varphi}, \tilde{\psi})\in\ \text{int}(X_+)$, there is $\rho>0$ such that $\tilde{\varphi} = \rho \varphi$ and $\tilde{\psi} = \rho\mu^{\frac{1}{\beta_2}} \psi$ in $\Omega$.
	
	Existence, nonexistence and uniqueness of nontrivial solutions to the system (\ref{1}) have been widely investigated during the last decades for $p=q=2$, $\alpha_1=\alpha_2=0$ and $s_1=s_2=1$, we refer for instance to \cite{CFM1992, FM1998, Fi, FF1994, HV1993,Mi1993,SZ1998}, where in particular notions of sub-superlinearity, sub-supercriticality and criticality have been introduced. Still in this context, the eigenvalue problem, i.e., $\beta_1 \beta_2=1$, was completely studied in \cite{Mon2000}. Recently some of these works have been extended to the fractional case, that is, for system of the type (\ref{1}) with $p=q=2$, $\alpha_1=\alpha_2=0$ and $s_1,s_2\in (0,1)$, existence results of positive solutions have been established when $\beta_1\beta_2>1$ in \cite{EM} for $s_1 \neq s_2$ and in \cite{EM1} for $s_1 = s_2$. The latter also proves existence and uniqueness of positive solution in the case that $\beta_1\beta_2<1$. Finally, when $\beta_1\beta_2 = 1$, the eigenvalue problem has been investigated in \cite{LM2020-2}.
	
	In the $(p,q)$-Laplacian case (when $s_1=s_2=1$) the system (\ref{1}) was studied in \cite{KR}, in that the authors has proved the existence of principal eigenvalues, more precisely, they proved a smooth curve of pairs $(\lambda, \mu)$ in $(0, \infty) \times(0, \infty)$ such that the quasilinear elliptic system possesses a solution pair $(u, v)$ consisting of nontrivial, nonnegative functions $u \in W_0^{1, p}(\Omega)$ and $v \in W_0^{1, q}(\Omega)$. In particular, the superlinear case when $\alpha_1=\alpha_2=0$ and $\beta_1 \beta_2>(p-1)(q-1)$ was treated in Ph. Clément, R. F. Manásevich and E. Mitidieri \cite{CMM_1993}. 
	
	There exists connection between principal eigenvalues and maximum principles as been investigated in \cite{AFL2009,Am2005,AL-G2019,CL-G2008,L-GM-M1994} for cooperative systems and in \cite{Sw1992} for non-cooperative systems and more recently in \cite{LM2020}, where system (\ref{1}) was analyzed in the special case when $p=q=2$, $\alpha_1=\alpha_2=0$ and $s_1=s_2=1$, however, instead of $\Delta$, a general second order elliptic operator was considered. For the nonlocal context, we refer \cite{LM2020-2}. In \cite{Leite2021} was extended the results of \cite{LM2020} for $p, q>1$, that is, was established the connection between principal spectral curves for systems (\ref{1})  and maximum and comparison principles related when $\alpha_1=\alpha_2=0$ and $\beta_1 \beta_2=(p-1)(q-1)$. Explicit lower bounds for principal eigenvalues for the system in terms of the measure of $\Omega$ were also proved.

In this paper we show the existence of principal eigenvalues of \eqref{1} and some of their qualitative properties. Namely, we prove that the set formed by these principal eigenvalues is given by the following smooth curve
\[
\mathcal{C}_1:=\left\{(\lambda,\mu)\in (\mathbb{R}_+^*)^2;\lambda^{\frac{1}{\sqrt{\beta_1(p-1-\alpha_1)}}}\mu^{\frac{1}{\sqrt{\beta_2(q-1-\alpha_2)}}}=\Lambda_0 \right\},
\]
for some $\Lambda_0>0$, which satisfies:

\begin{enumerate}
			\item[(P1)] $(\lambda, \mu) \in \mathbb{R}_+\times\mathbb{R}_+$ is a principal eigenvalue of the system (\ref{1}) if, and only if, $(\lambda,\mu)\in \mathcal{C}_1$;
			
			\item[(P2)] The curve $\mathcal{C}_1$ is simple in $\text{int}(X_+)$, that is, $(\lambda,\mu)$ is simple in $\text{int}(X_+)$ for all $(\lambda,\mu)\in \mathcal{C}_1$;
			\item[(P3)] Suppose that $\alpha_1 = \alpha_2 = 0$ and let $(u,v) \in X$ be an eigenfunction associated to $(\lambda,\mu) \in \mathcal{C}_1$. Then either $(u,v) \in \text{int}(X_+)$ or $(-u,-v) \in\text{int}(X_+)$.
		\end{enumerate}

More precisely we will prove the following result in the homogeneous case, that is, when $\alpha_1 + \beta_1 = p-1$ e $\alpha_2+\beta_2 =q-1$.
	\begin{teor}\label{Th01}
		Suppose that
		\begin{equation}\label{iden2}
			\alpha_1 + \beta_1 = p-1 \ \ \text{and} \ \ \alpha_2 + \beta_2 = q-1.
		\end{equation}
		Then there exists $\Lambda>0$ and a couple $(u_1,v_1) \in \text{int}(X_+)$ such that the system \eqref{1} admits a positive solutions $(u,v) \in X_+$ corresponding to some $(\lambda, \mu) \in (\mathbb{R}_+^{*})^2$ if, and only if, 
		
		\[
		\lambda^{\frac{1}{\beta_1}} \mu^{\frac{1}{\beta_2}} = \Lambda.
		\]
		Furthermore, $u = cu_1$ and $v=cv_1$ with $c>0$.
	\end{teor}

	In the general case, we will prove the following results:
	
\begin{teor}\label{teor2}
		Let $\alpha_1 < p-1$, $\alpha_2 < q-1$ and 
		\begin{equation}\label{cond1}
			\beta_1\beta_2 = (p-1-\alpha_1)(q-1-\alpha_2).
		\end{equation}
		
		Then there exists $\Lambda_0 >0$ and a couple $(u_0,v_0) \in \text{int}(X_+)$ such that system \eqref{1} admits a positive solution $(u,v) \in X_+$ corresponding to some $(\lambda, \mu) \in (\mathbb{R}_+^{*})^2$ if, and only if, 
		\begin{equation}\label{cond2}
			\lambda^{\frac{1}{\sqrt{\beta_1(p-1-\alpha_1)}}} \mu^{\frac{1}{\sqrt{\beta_2(q-1-\alpha_2)}}} = \Lambda_0.
		\end{equation}
		
		Furthermore, there is $\rho >0$ such that $u=\rho u_0$ and $v = \rho \mu^{\frac{1}{\beta_2}}v_0$.
	\end{teor}

	\begin{teor}\label{teor3}
		Suppose that $\alpha_1<p-1$, $\alpha_2<q-1$ and \eqref{cond1} hold. Then properties (P1), (P2) and (P3) are satisfied.
	\end{teor}
	
	Another important property about the principal eigenvalues concerns the lower estimate. Such estimates have been widely investigated in recent decades (see for example in \cite{BeNiVa, CaLo,CL-G2008, Pinasco, LM2020, L, L1, L2, Lo, Protter}). Here, we will use the $L^{\infty}$-bound (\ref{abp}) to obtain a lower estimate associated with the system \eqref{1}. Namely:

\begin{teor} \label{lower} Let $a,b\in L^\infty(\Omega)$ and $\mathcal{C}_1$ be the principal curve associated to \eqref{1}. Suppose that $\alpha_1<p-1$, $\alpha_2<q-1$ and \eqref{cond1} hold. Then 

\begin{eqnarray}
&& \Lambda_0 \geq \frac{1}{C_{s_1,p}^{\frac{p-1}{\theta}}C_{s_2,q}^{\frac{q-1}{\zeta}}d^{\frac{ps_1}{\theta}+\frac{qs_2}{\zeta}} \|a\|^{\frac{1}{\theta}}_{L^\infty(\Omega)} \|b\|^{\frac{1}{\zeta}}_{L^\infty(\Omega)}},  \label{lb1} 
\end{eqnarray}
where $d=\operatorname{diam }(\Omega)$, $\theta:=\sqrt{\beta_1(p-1-\alpha_1)}$, $\zeta:=\sqrt{\beta_2(q-1-\alpha_2)}$. In particular,
\[
\lim_{d\downarrow 0}\Lambda_0=+\infty.
\]
\end{teor}

Notice that we get an explicit lower estimate of $\Lambda_0$ in terms of the diameter of $\Omega$, weighted functions $a, b \in L^\infty(\Omega)$ and the explicit constants $C_{s_1,p}, C_{s_2,q}$.

	The second part of the work is devoted to the study of maximum and comparison principles related to the following system:
\begin{equation}\label{1.3}
\left\{
\begin{array}{llll}
(-\Delta_p)^{s_1}u = \lambda a(x)\vert v\vert^{\beta_1-1}v & {\rm in} \ \ \Omega,\\
(-\Delta_q)^{s_2}v = \mu b(x)\vert u\vert^{\beta_2-1}u & {\rm in} \ \ \Omega,\\
u= v=0 & {\rm in} \ \ \mathbb{R}^N\setminus\Omega,
\end{array}
\right.
\end{equation}
where $\beta_1, \beta_2 > 0$ with $\beta_1 \beta_2 = (p-1)(q-1)$, and $(\lambda, \mu) \in \mathbb{R}^2$. 

\begin{rem}
Notice that the system (\ref{1.3}) is a particular case of the system (\ref{1}) with $\alpha_1=\alpha_2=0$.
\end{rem}

We define
\[
\widetilde{W}^{s,m}(\Omega):=\left\{ u \in L^m_{loc}(\mathbb{R}^N);\ \int\limits_{\Omega}\int\limits_{\mathbb{R}^N}\frac{|u(x) - u(y)|^m}{|x-y|^{N + sm}}\ dxdy<\infty\right\}.
\]

The space $\mathcal{W}^{s,m}(\Omega)$ is the space of all function $u \in L^m_{loc}(\mathbb{R}^N)$ such that for
any bounded $\Omega'\subseteq\Omega$ there is an open set $U\supset\supset\Omega'$ so that $u\in W^{s,m}(U)$, and
\[
\int\limits_{\mathbb{R}^N}\frac{\vert u(x)\vert^{m-1}}{(1+\vert x\vert)^{N+sm}}\ dx<+\infty,
\]
in particular 
\[
W^{s,m}_0(\Omega)\subset\mathcal{W}^{s,m}(\Omega)\cap\widetilde{W}^{s,m}(\Omega),
\]
for more details see \cite{Iannizzotto}.

A couple $(u,v)\in \widetilde{W}^{s_1,p}(\Omega)\times \widetilde{W}^{s_2,q}(\Omega)$ is a weak solution of the system 
\begin{equation}\label{supersol}
\left\{
\begin{array}{llll}
(-\Delta_p)^{s_1} u \geq \lambda a(x)\vert v\vert^{\beta_1-1}v & {\rm in} \ \ \Omega,\\
(-\Delta_q)^{s_2} v \geq \mu b(x)\vert u\vert^{\beta_2-1}u & {\rm in} \ \ \Omega,\\
u, v\geq 0 & {\rm in} \ \ \mathbb{R}^N\setminus\Omega,
\end{array}
\right.
\end{equation}
if the inequalities
\[
C(N,s_1,p)\int\limits_{\mathbb{R}^N}\int\limits_{\mathbb{R}^N}\frac{\vert u(x)-u(y)\vert^{p-2}(u(x)-u(y))(f_1(x)-f_1(y))}{\vert x-y\vert^{N+s_1p}}\ dxdy\geq\int\limits_{\Omega}\lambda a(x)\vert v\vert^{\beta_1-1}vf_1\ dx
\]
and 
\[
C(N,s_2,q)\int\limits_{\mathbb{R}^N}\int\limits_{\mathbb{R}^N}\frac{\vert v(x)-v(y)\vert^{q-2}(v(x)-v(y))(f_2(x)-f_2(y))}{\vert x-y\vert^{N+s_2q}}\ dxdy\geq\int\limits_{\Omega}\mu b(x)\vert u\vert^{\beta_2-1}uf_2\ dx
\]
hold for every $0\leq f_1\in W^{s_1,p}_0(\Omega)$, $0\leq f_2\in W^{s_2,q}_0(\Omega)$ and $u, v\geq 0$ a.e. in $\mathbb{R}^N\setminus\Omega$.

By weak maximum principle, denoted by {\bf (WMP)}, we mean that for any weak solution $(u,v)\in \mathcal{W}^{s_1,p}(\Omega)\cap\widetilde{W}^{s_1,p}(\Omega)\cap C(\mathbb{R}^N)\times \mathcal{W}^{s_2,q}(\Omega)\cap\widetilde{W}^{s_2,q}(\Omega)\cap C(\mathbb{R}^N)$ of the system (\ref{supersol}) with $(u,v)$ satisfying 
\begin{equation}\label{h1}
\int\limits_{\Omega}\int\limits_{\Omega}\frac{\left\vert\vert\varphi(x)-\varphi(y)\vert^{p-2}(\varphi(x)-\varphi(y))-\vert u(x)-u(y)\vert^{p-2}(u(y)-u(x))\right\vert}{\vert x-y\vert^{N+s_1p}}\ dxdy<\infty
\end{equation}
and 
\begin{equation}\label{h2}
\int\limits_{\Omega}\int\limits_{\Omega}\frac{\left\vert\vert\psi(x)-\psi(y)\vert^{q-2}(\psi(x)-\psi(y))-\vert v(x)-v(y)\vert^{q-2}(v(y)-v(x))\right\vert}{\vert x-y\vert^{N+s_2q}}\ dxdy<\infty
\end{equation}
for all pair $(\varphi,\psi)$ of positive eigenfunction associated to $(\lambda_1, \mu_1)$, where 
\[
		\left\{
		\begin{array}{llll}
			(\lambda_1, \mu_1)\in\mathcal{C}_1\ {\rm with}\ \frac{\mu}{\lambda}=\frac{\mu_1}{\lambda_1} & {\rm when} \ \lambda\ {\rm and}\ \mu\ {\rm have\ the\ same\ sign},\\
			(\lambda_1, \mu_1)\in\mathcal{C}_1\ {\rm with}\ \frac{\mu}{\lambda}=-\frac{\mu_1}{\lambda_1} & {\rm when} \ \lambda\ {\rm and}\ \mu\ {\rm have\ the\ diferent\ sign},
		\end{array}
		\right.
	\]
we have $u, v \geq  0$ in $\Omega$. Moreover, if at least, $u$ or $v$ is positive in $\Omega$, we say that the strong maximum principle, denoted by {\bf (SMP)}, corresponding to \eqref{1.3} holds in $\Omega$. Now, if $\lambda, \mu > 0$, {\bf (SMP)} can be rephrased as $u, v > 0$ in $\Omega$. 

The following result characterize completely in terms of the curve $\mathcal{C}_1$ the set of $(\lambda, \mu) \in \mathbb{R}^2$ such that {\bf (WMP)} and {\bf (SMP)} hold in $\Omega$. Precisely:
\begin{figure}[ht]
\centering
\includegraphics[scale=1]{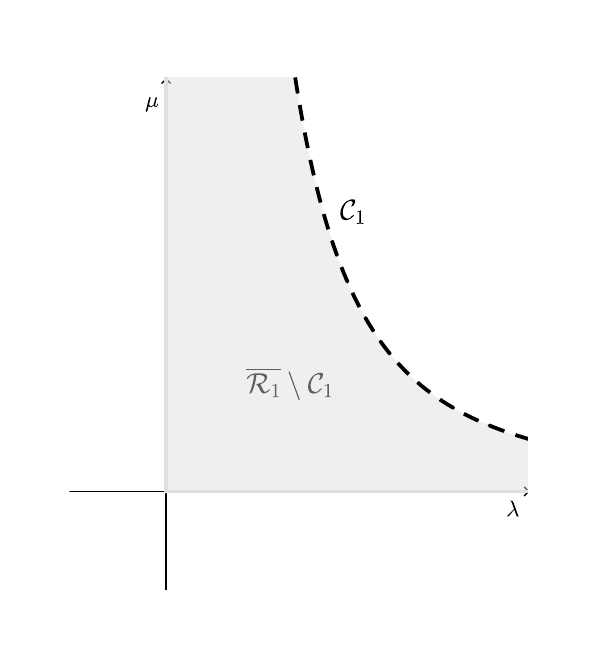}
\caption{Region of validity for {\bf (WMP)} and {\bf (SMP)}.}
\end{figure}

\begin{teor} \label{MP}
Let $(\lambda, \mu) \in \mathbb{R}^2$ and $\mathcal{R}_1$ be the open region in the first quadrant below $\mathcal{C}_1$. The following assertions are equivalent:

\begin{itemize}
\item[{\rm (i)}] $(\lambda, \mu) \in \overline{\mathcal{R}_1} \setminus \mathcal{C}_1$;
\item[{\rm (ii)}] {\bf (WMP)} associated to \eqref{1.3} holds in $\Omega$;
\item[{\rm (iii)}] {\bf (SMP)} associated to \eqref{1.3} holds in $\Omega$.
\end{itemize}
\end{teor}

The pairs $(u,v)$ and $(z,w)$ in $\widetilde{W}^{s_1,p}(\Omega)\times \widetilde{W}^{s_2,q}(\Omega)$ are weak solutions of the system 
\begin{equation}\label{comp}
\left\{
\begin{array}{llll}
0\leq (-\Delta_p)^{s_1} u  - \lambda a(x)\vert v\vert^{\beta_1-1}v \leq (-\Delta_p)^{s_1} z  - \lambda a(x)\vert w\vert^{\beta_1-1}w & {\rm in} \ \ \Omega,\\
0\leq (-\Delta_q)^{s_2} v  - \mu b(x)\vert u\vert^{\beta_2-1}u \leq (-\Delta_q)^{s_2} w  - \mu b(x)\vert z\vert^{\beta_2-1}z & {\rm in} \ \ \Omega,\\
u \leq z\ \ {\rm and}\ \ v \leq  w & {\rm in} \ \ \mathbb{R}^N\setminus\Omega,
\end{array}
\right.
\end{equation}
if the inequalities
\begin{eqnarray*}
0&\leq &C(N,s_1,p)\int\limits_{\mathbb{R}^N}\int\limits_{\mathbb{R}^N}\frac{\vert u(x)-u(y)\vert^{p-2}(u(x)-u(y))(f_1(x)-f_1(y))}{\vert x-y\vert^{N+s_1p}}\ dxdy-\int\limits_{\Omega}\lambda a(x)\vert v\vert^{\beta_1-1}vf_1\ dx\\
&\leq &C(N,s_1,p)\int\limits_{\mathbb{R}^N}\int\limits_{\mathbb{R}^N}\frac{\vert z(x)-z(y)\vert^{p-2}(z(x)-z(y))(f_1(x)-f_1(y))}{\vert x-y\vert^{N+s_1p}}\ dxdy-\int\limits_{\Omega}\lambda a(x)\vert w\vert^{\beta_1-1}wf_1\ dx
\end{eqnarray*}
and 
\begin{eqnarray*}
0&\leq &C(N,s_2,q)\int\limits_{\mathbb{R}^N}\int\limits_{\mathbb{R}^N}\frac{\vert v(x)-v(y)\vert^{q-2}(v(x)-v(y))(f_2(x)-f_2(y))}{\vert x-y\vert^{N+s_2q}}\ dxdy-\int\limits_{\Omega}\mu b(x)\vert u\vert^{\beta_2-1}uf_2\ dx\\
&\leq &C(N,s_2,q)\int\limits_{\mathbb{R}^N}\int\limits_{\mathbb{R}^N}\frac{\vert w(x)-w(y)\vert^{q-2}(w(x)-w(y))(f_2(x)-f_2(y))}{\vert x-y\vert^{N+s_2q}}\ dxdy-\int\limits_{\Omega}\mu b(x)\vert z\vert^{\beta_2-1}zf_2\ dx
\end{eqnarray*}
hold for every $0\leq f_1\in W^{s_1,p}_0(\Omega)$, $0\leq f_2\in W^{s_2,q}_0(\Omega)$ and $u\leq z$, $v\leq w$ a.e. in $\mathbb{R}^N\setminus\Omega$.

Note that weak and strong comparison principles, denoted respectively by {\bf (WCP)} and {\bf (SCP)}, are essential tools for establish the positivity and uniqueness of solutions for elliptic problems, among others, to certain counterparts of (\ref{1.3}). We say that {\bf (WCP)} holds in $\Omega$ if, for any weak solutions $(u,v)$ and $(z,w)$ in $\mathcal{W}^{s_1,p}(\Omega)\cap\widetilde{W}^{s_1,p}(\Omega)\cap C(\mathbb{R}^N)\times \mathcal{W}^{s_2,q}(\Omega)\cap\widetilde{W}^{s_2,q}(\Omega)\cap C(\mathbb{R}^N)$ of the system (\ref{comp}) with the pairs $(u,v)$ and $(z,w)$ satisfying (\ref{h1}) and (\ref{h2}),
\begin{equation}\label{h3}
\int\limits_{\Omega}\int\limits_{\Omega}\frac{\left\vert\vert z(x)-z(y)\vert^{p-2}(z(x)-z(y))-\vert \delta u(x)-\delta u(y)\vert^{p-2}(\delta u(x)-\delta u(y))\right\vert}{\vert x-y\vert^{N+s_1p}}\ dxdy<\infty
\end{equation}
and 
\begin{equation}\label{h4}
\int\limits_{\Omega}\int\limits_{\Omega}\frac{\left\vert\vert w(x)-w(y)\vert^{q-2}(w(x)-w(y))-\vert \delta^{\omega}v(x)-\delta^{\omega}v(y)\vert^{q-2}(\delta^{\omega}v(x)-\delta^{\omega}v(y))\right\vert}{\vert x-y\vert^{N+s_2q}}\ dxdy<\infty
\end{equation}
for all $0<\delta\leq 1$, where $\omega:=\frac{p-1}{\beta_1}$, one has $u \leq z$ and $v \leq w$ in $\Omega$. Besides, the strong comparison principle {\bf (SCP)} is said to hold in $\Omega$ if either $u \equiv z$ and $v \equiv w$ in $\Omega$ or at least $u < z$ in $\Omega$ or $v < w$ in $\Omega$. In the case that $\lambda, \mu >0$, {\bf (SCP)} in $\Omega$ can be rephrased as either $u \equiv z$ and $v \equiv w$ in $\Omega$ or $u < z$ and $v < w$ in $\Omega$.

Our next theorem classify completely the {\bf (WCP)} and {\bf (SCP)} associated to (\ref{1.3}) in terms of the principal curve $\mathcal{C}_1$. Namely:

\begin{teor} \label{CP} Let $\mathcal{R}_1$ be as in Theorem \ref{MP} and $(\lambda, \mu) \in \mathbb{R}^2$. The following assertions are equivalent:

\begin{itemize}
\item[{\rm (i)}] $(\lambda, \mu) \in \overline{\mathcal{R}_1} \setminus \mathcal{C}_1$;
\item[{\rm (ii)}] {\bf (WCP)} associated to \eqref{1.3} holds in $\Omega$;
\item[{\rm (iii)}] {\bf (SCP)} associated to \eqref{1.3} holds in $\Omega$.
\end{itemize}
\end{teor}

The validity of strong comparison principle for problems involving the fractional $m$-Laplacian operator is delicate (see \cite{CFGQ, Mosconi}). For example,  the operator $(-\Delta_m)^{s}$ satisfies the strong comparison principle only under the additional condition zero Dirichlet boundary values in $\mathbb{R}^N\setminus\Omega$ and positives in $\Omega$ (see \cite{Mosconi}) or under the condition (\ref{h3}) with $\delta=1$ (see \cite{CFGQ}). In our work, we use the strong comparison principle obtained in \cite{Mosconi} to prove the existence of a spectral principal curve corresponding to the system (\ref{1}). On the other hand, to obtain the maximum and comparison principles associated to (\ref{1.3}), it was possible to apply only the strong comparison principle obtained in \cite{CFGQ}. So, similar conditions of the scalar context are necessary for the systems context. 

A topic of great interest theme within the elliptic PDEs theory concerns maximum principles depending on domain $\Omega$. Various of the contributions can be found in the list of references \cite{BeNiVa, CaLo, LM2020, L1, L2, Lo}.

Inspired by these ideas and by using the lower estimate (\ref{lb1}), we characterize maximum principles corresponding to \eqref{1.3} depending on diameter of domain. Namely:

\begin{teor} \label{sm} Let $\theta$ and $\zeta$ be as in Theorem \ref{lower}, with $\alpha_1=0=\alpha_2$. Define
\[
\eta := \frac{1}{\left[\lambda^{\frac{1}{\theta}}\mu^{\frac{1}{\zeta}}C_{s_1,p}^{\frac{p-1}{\theta}}C_{s_2,q}^{\frac{q-1}{\zeta}} \|a\|^{\frac{1}{\theta}}_{L^\infty(\Omega)} \|b\|^{\frac{1}{\zeta}}_{L^\infty(\Omega)}\right]^{\frac{\beta_1\beta_2}{ps_1\theta+qs_2\zeta}}}.
\]
The following assertions are equivalent:

\begin{itemize}
\item[{\rm (i)}] $\lambda \geq 0$ and $\mu \geq 0$;
\item[{\rm (ii)}] {\bf (WMP)} corresponding to \eqref{1.3} holds in $\Omega$ provided that $d < \eta$;
\item[{\rm (iii)}] {\bf (SMP)} corresponding to \eqref{1.3} holds in $\Omega$ provided that $d < \eta$.
\end{itemize}
\end{teor}

In order to make the paper more organized, we show each theorem in its corresponding section. 
	
\section{Proof of Theorem \ref{Th01}}
\mbox{}

We consider the cone $(C^0_{s}(\overline{\Omega}))_+ :=\bigg\{f \in C^0_s(\overline{\Omega}) ; \ \ f \geq 0 \ \ \text{in} \ \ \Omega\bigg\}$ which is nonempty and has topological interior 
\[
\text{int} \left[\left(C^{0}_{s}(\overline{\Omega})\right)_+\right] = \bigg\{ v \in C^0_s(\overline{\Omega}); \ \ v > 0 \ \ \text{in} \ \ \Omega \ \ \text{and} \ \ \inf\frac{v}{d_{\Omega}^s} > 0 \bigg\},
\]
where $d_{\Omega}$ is the distance function. 
		
We consider also $X = C^{0}_{s_1}(\overline{\Omega}) \times C^{0}_{s_2}(\overline{\Omega})$, which is a strongly
ordered Banach space endowed with the natural norm and ordering for pairs of functions $(f,g) \in X$. Its positive cone 
$$X_+ = \bigg\{(f,g) \in X; \ \ f \geq 0 \ \ \text{and} \ \ g \geq 0 \ \ \text{in} \ \ \Omega \bigg\}$$
is normal and has nonempty topological interior 
\[
\text{int}(X_+) = \text{int}\left[\left(C^0_{s_1}(\overline{\Omega})\right)_+\right] \times \text{int}\left[\left(C^0_{s_2}(\overline{\Omega})\right)_+\right].
\]

Finally define the map $S: X \longrightarrow X$ by $S(u,v)=(\tilde{u},\tilde{v})$ with
\[
\tilde{u} = T_p^{s_1}(a(x)\vert u \vert^{\alpha_1}\vert v \vert^{\beta_1-1}v) \ \ \text{and} \ \  \tilde{v} = T_q^{s_2}(b(x)\vert v \vert^{\alpha_2}\vert u \vert^{\beta_2-1}u)
\]
for $(u,v) \in X$.
		
The lemma below describes the interaction between principal eigenvalues and the different homogeneities of $S_1$ and $S_2$, the two components of $S$.
\begin{lema}\label{Lema1}
\begin{enumerate}
\item[(i)] A couple $(u,v) \in X_+\backslash \{0\}$ is a weak solution of the system \eqref{1} for some $(\lambda, \mu) \in (\mathbb{R}^{*}_+)^2$ if, and only if, $(u,v) \in \text{int}(X_+)$ and $S(u,v) = ( \lambda^{\frac{-1}{p-1}}u, \mu^{\frac{-1}{q-1}}v)$;
			
\item[(ii)] For all $\rho, \sigma \in \mathbb{R}_+$ and for all $(u,v) \in X_+$ we have 
\[
S(\rho u,\sigma v) = \big((\rho^{\alpha_1}\sigma^{\beta_1})^{\frac{1}{p-1}}S_1(u,v),(\rho^{\alpha_2}\sigma^{\beta_2})^{\frac{1}{q-1}}S_2(u,v)\big);
\]
			
\item[(iii)] If $(u,v) \in X_+$ solves \eqref{1} with $(\lambda_0,\mu_0)$ in place of $(\lambda,\mu)$, then for any $\rho,\sigma >0$, the pair $(\rho u,\sigma v)$ solves \eqref{1} with 
\[
\lambda = \rho^{p-1-\alpha_1}\sigma^{-\beta_1}\lambda_0 \ \ \text{and} \ \ \mu= \rho^{q-1-\alpha_2}\sigma^{-\beta_2}\mu_0.
\]
\end{enumerate}
\end{lema}  
	\n {\bf Proof.}
	Notice that, give $t > 0$ and $m>1$, we have
	\begin{equation}\label{iden1}
		(-\Delta_m)^s(tu) = t^{m-1} (-\Delta_m)^{s} u, \ \ \text{for all} \ \ u \in W^{s,m}(\Omega).
	\end{equation}

	$(i)$ Note that, if $(u,v) \in X_+ \backslash \{0\}$ is a weak solution of \eqref{1} then, by Hopf's Lemma and maximum principle, see \cite{Ian} and \cite{Pezzo}, we have 
	\[
	u,v > 0 \ \ \text{in} \ \ \Omega \ \ \text{and} \ \ \frac{u(x)}{d_{\Omega}^{s_1}(x)}, \frac{v(x)}{d_{\Omega}^{s_2}(x)} > 0 \ \ \text{for all} \ \ x \in \mathbb{R}^N
	\]
	and thus, $(u,v) \in \text{int}(X_+)$.
	
	Moreover, by identity \eqref{iden1} results
	\[
	(-\Delta_p)^{s_1}u = \lambda a(x) \vert u \vert^{\alpha_1} \vert v \vert^{\beta_1 -1}v  \ \ \text{and} \ \ (-\Delta_q)^{s_2}v = \mu b(x) \vert v\vert^{\alpha_2} \vert u \vert^{\beta_2 -1}u 
	\]
	if, and only if,
	\[
	(-\Delta_p)^{s_1}(\lambda^{\frac{-1}{p-1}}u) = a(x) \vert u \vert^{\alpha_1} \vert v \vert^{\beta_1 -1}v \ \ \text{and} \ \ (-\Delta_q)^{s_2}(\lambda^{\frac{-1}{q-1}}v) = a(x) \vert v \vert^{\alpha_2} \vert u \vert^{\beta_2 -1}u.
	\]
	Therefore, $(u,v) \in X_+\backslash\{0\}$ is a solution of \eqref{1} if, and only if, $S(u,v) = \big(\lambda^{\frac{-1}{p-1}}u, \mu^{\frac{-1}{q-1}}v\big)$.
	
	$(ii)$ Let $\rho , \sigma \in \mathbb{R}^+$ and $(u,v) \in X_+$. Denoting $\tilde{u} = T_p^{s_1}(a(x) \vert \rho u \vert^{\alpha_1} \vert \sigma v \vert^{\beta_1 -1}\sigma v)$ we have
	
	\begin{align*}
		(-\Delta_p)^{s_1}\tilde{u} = a(x) \vert \rho u \vert^{\alpha_1} \vert \sigma v \vert^{\beta_1 -1}\sigma v
	\end{align*}
	and thus
	\[
	(-\Delta_p)^{s_1}((\rho^{-\alpha_1}\sigma^{-\beta_1})^{\frac{1}{p-1}}\tilde{u}) = a(x) \vert u \vert^{\alpha_1} \vert v \vert^{\beta_1 -1} v,
	\]
	that is,
	\[
	(\rho^{-\alpha_1}\sigma^{-\beta_1})^{\frac{1}{p-1}}\tilde{u} = T_p^{s_1}(a(x)\vert u \vert^{\alpha_1} \vert v \vert^{\beta_1-1}v) = S_1(u,v).
	\]
	Therefore
	\[
	\tilde{u} = (\rho^{\alpha_1}\sigma^{\beta_1})^{\frac{1}{p-1}}S_1(u,v).
	\]
	
	Analogously, show that, if $\tilde{v} = T_q^{s_2}(b(x) \vert \sigma v \vert^{\alpha_2} \vert \rho u \vert^{\beta_2 -1}\rho u)$ then
	\[
	\tilde{v} = (\rho^{\beta_2}\sigma^{\alpha_2})^{\frac{1}{q-1}}S_2(u,v).
	\]
	
	$(iii)$ For this item just calculate directly observing that, if
	\[
	(-\Delta_p)^{s_1}u = \lambda_0 a(x) \vert u \vert^{\alpha_1} \vert v \vert^{\beta_1 -1} v
	\]
	then,
	\begin{align*}
		(-\Delta_p)^{s_1}(\rho u) &= \rho^{p-1}(-\Delta_p)^su \\
		&= \lambda_0 \rho^{p-1}a(x) \vert u \vert^{\alpha_1} \vert v \vert^{\beta_1 -1} v \\
		&= \lambda_0 \rho^{p-1- \alpha_1} \sigma^{-\beta_1}a(x) \vert \rho u \vert^{\alpha_1} \vert \sigma v \vert^{\beta_1 -1} \sigma v
	\end{align*}
	and in the same way we find for $\mu = \rho^{-\beta_2}\sigma^{q-1-\alpha_2}\mu_0$.
	\qed
	
{\bf Proof of Theorem \ref{Th01}:}\\
\mbox{}

 Using the item $(i)$ of the Lemma \ref{Lema1}, we have for all $t \geq 0$,
	\[
	S(tu,tv) = \bigg((t^{\alpha_1}t^{\beta_1})^{\frac{1}{p-1}}S_1(u,v), (t^{\alpha_2}t^{\beta_2})^{\frac{1}{q-1}}S_2(u,v)\bigg) = \bigg( t^{\frac{\alpha_1 + \beta_1}{p-1}} S_1(u,v), t^{\frac{\alpha_2 + \beta_2}{q-1}}S_2(u,v)\bigg)
	\]
	and thus, $S(tu,tv) = tS(u,v)$) if, and only if, $\alpha_1 + \beta_1 = p-1$ and $\alpha_2 + \beta_2 = q-1$. Therefore, $S$ is homogeneous if, and only if, \eqref{iden2} is satisfied.
	
	Consider $(u_i,v_i) \in X$, with $i=1,2$. If $0 \leq u_1\leq u_2$ and $0 \leq v_1 \leq v_2$ in $\Omega$, with $u_1 \not\equiv u_2$ and $v_1 \not\equiv v_2$ then, 
	\[
	a(x) \vert u_1\vert^{\alpha_1} \vert v_1 \vert^{\beta_1-1}v_1 \leq a(x) \vert u_2\vert^{\alpha_1} \vert v_2 \vert^{\beta_1-1}v_2
	\]
	and thus if $\tilde{u}_i = T_p^{s_1}(a(x) \vert u_i\vert^{\alpha_1} \vert v_i \vert^{\beta_1-1}v_i)$ we have 
	\[
	(-\Delta_p)^{s_1} \tilde{u}_1 \leq (-\Delta_p)^{s_1}\tilde{u}_2  \ \ in \ \ \Omega.
	\]
	By strong comparison principle, see Theorem 2.7 in \cite{Mosconi}, results that $\tilde{u}_1 < \tilde{u}_2$ in $\Omega$. 
	
	Analogously we show that $\tilde{v}_1 < \tilde{v}_2$ and thus results that $S(u_1,v_1) - S(u_2,v_2) \in   \text{int}(X_+)$.
	
	Finally, by a regularity result interior which can be found in (\cite{Iannizzotto},Theorem 1.1), it follows that the mapping $T_p^{s}: L^{\infty}(\Omega) \longrightarrow C^{\alpha}(\overline{\Omega})$ is continuous and bounded; that is, it maps bounded sets into bounded sets. Then, by Arzelá-Ascoli's theorem, $T_p^{s}: L^{\infty}(\Omega) \longrightarrow C^{\alpha'}(\overline{\Omega})$ is
	compact whenever $0 < \alpha' \leq \alpha$, that is, it maps bounded sets into sets with compact closure.
	By Krein-Rutman theorem, see (\cite{KR}, Theorem A.2),  there exists $\Lambda_1 \in \mathbb{R}$ and $ (u_1,v_1) \in \text{int}(X_+)$ such that,  $S(u_1,v_1) = \Lambda_1(u_1,v_1)$. From Lemma \ref{Lema1} results that the couple $(\lambda_0, \mu_0) = (\Lambda_1^{1-p}, \Lambda_1^{1-q})$ is a eigenvalue for system \eqref{1}. 
	
	Note that, using the conditions, $\alpha_1 + \beta_1 = p-1$ and $\alpha_2 + \beta_2 = q-1$ we have, for any $(\lambda, \mu)$ satisfying 
	\[
	\left(\frac{\lambda}{\Lambda_1^{1-p}}\right)^{\frac{1}{\beta_1}} \left(\frac{\mu}{\Lambda_1^{1-q}}\right)^{\frac{1}{\beta_2}} = 1
	\]
	the couple $\rho = \left(\frac{\lambda}{\Lambda_1^{1-p}}\right)^{\frac{1}{\beta_1}}$ and $\sigma=1$ satisfy the item $(iii)$ of the Lemma \ref{Lema1}.
	
	Reciprocally, if the system \eqref{1} possesses a positive weak solution $(u,v) \in \text{int}(X_+)$ corresponding to some $(\lambda, \mu)$, taking $\rho, \sigma >0$ such that $((\frac{\rho}{\sigma})^{\beta_1}\lambda)^{\frac{-1}{p-1}} = ((\frac{\sigma}{\rho})^{\beta_2}\mu)^{\frac{-1}{q-1}} = \Lambda_0$, the Lemma \ref{Lema1} assures us that 
	\begin{align*}
		S(\rho u, \sigma v) &= ((\rho^{\alpha_1}\sigma^{\beta_1})^{\frac{1}{p-1}}S_1(u,v),(\rho^{\beta_2}\sigma^{\alpha_2})^{\frac{1}{q-1}}S_2(u,v))\\
		&= ((\rho^{p-1-\beta_1}\sigma^{\beta_1})^{\frac{1}{p-1}}\lambda^{-\frac{1}{p-1}} u,(\rho^{\beta_2}\sigma^{q-1-\beta_2})^{\frac{1}{q-1}}\mu^{-\frac{1}{q-1}} v)\\
		&= \left(\left( \left(\frac{\rho}{\sigma}\right)^{\beta_1}\lambda\right)^{-\frac{1}{p-1}}\rho u, \left( \left(\frac{\sigma}{\rho}\right)^{\beta_1}\mu\right)^{-\frac{1}{q-1}}\sigma v\right)\\
		&= \Lambda_0\left(\rho u, \sigma v\right).
	\end{align*}

	By Krein-Rutman theorem, see (\cite{KR}, Theorem A.2), results that $\Lambda_0 = \Lambda_1$. Thus, talking 
	\[
	\Lambda = \Lambda_1^{\frac{-(p-1)\beta_2}{\beta_1\beta_2 q-1}}
	\]
	we complete the proof of the theorem.
	\qed

	\section{Proof of Theorem \ref{teor2}}
\mbox{}	 

	Let us assume that $\alpha_1 < p-1$ and $\alpha_2 < q-1$
and consider the new mapping $T: X_+ \longrightarrow X_+$ defined by, 
	\[
	T(u,v) = (J_1(v), J_2(u)), \ \ \text{for all} \ \ (u,v) \in X_+
	\]
	where $J_1(v)$ is the unique weak solution of the problem
	\[
		\left\{
		\begin{array}{rrll}
			(-\Delta_p)^{s_1} u &=& a(x) \vert u\vert^{\alpha_1}v^{\beta_1} & {\rm in} \ \ \Omega,\\
			u &=& 0  & {\rm in} \ \ \mathbb{R}^{N}\setminus\Omega
		\end{array}
		\right.
	\]
	and $J_2(u)$ is the unique weak solution of the problem
	\[
		\left\{
		\begin{array}{rrll}
			(-\Delta_q)^{s_2} u &=& b(x) \vert v\vert^{\alpha_2}u^{\beta_2} & {\rm in} \ \ \Omega,\\
			v &=& 0  & {\rm in} \ \ \mathbb{R}^{N}\setminus\Omega.
		\end{array}
		\right.
	\]
	
	The existence of solutions follows the classical minimization methods, more specifically the Lagrange multipliers method. The uniqueness follows from a convexity argument. Using regularity and strong maximum principle results, see \cite{Pezzo}, we conclude that $(J_1(v),J_2(u)) \in \text{int}(X_+)$.
	
	The idea in this case is to apply the Krein-Rutman theorem to the components of the mapping $T^2$. Notice that,
	
	\[
	T^2(u,v) = T(J_1(v), J_2(u)) = (J_1 \circ J_2(u), J_2 \circ J_1(v)).
	\]
	
	\begin{lema}\label{Lema2}
		The mapping $J_i: \left(C^0_{s_i}(\Omega)\right)_+ \longrightarrow \left(C^0_{s_i}(\Omega)\right)_+$ is nondecreasing for $i=1,2$.
	\end{lema}
	\n {\bf Proof.} 
	We consider $J_1$. Let $v_1,v_2 \in  C^0_{s_1}(\Omega)_+$, with $0 \leq v_1 \leq v_2$ and $v_1 \not\equiv 0$. Denote by $m_i(x) = a(x) v_i^{\beta_1}$ and $u_i = J_1(v_i)$. Since $v_1 \leq v_2$ we have $m_1(x) \leq m_2(x)$ and thus,
	\begin{equation*}
		\left\{
		\begin{array}{rrll}
		(-\Delta_p)^{s_1} u_2 &=& a(x) \vert u_2\vert^{\alpha_1} v_2^{\beta_1} = m_2(x) \vert u_2\vert^{\alpha_1} \geq m_1(x) \vert u_2 \vert^{\alpha_1} & {\rm in} \ \Omega,\\
		u_2&=& 0  & {\rm in} \ \ \mathbb{R}^{N}\setminus\Omega
		\end{array}
		\right.
	\end{equation*}
	in other words $u_2$ is a upper solutions of the problem 
	\begin{equation}\label{P2}
	\left\{
	\begin{array}{rrll}
	(-\Delta_p)^{s_1} u &=& m_1(x) \vert u \vert^{\alpha_1} \ &{\rm in} \ \Omega,\\   u_2 &=& 0 \ &{\rm in} \ \mathbb{R}^{N}\setminus\Omega.
\end{array}
\right.
\end{equation}
	
	Let us denote by $\varphi$ the positive eigenfunction of the Dirichlet fractional $p$-Laplacian with weight $m_1$, that is, there exists $\lambda_1(p) > 0$ such that 
	\[
	\left\{
	\begin{array}{rrll}
	(-\Delta_p)^{s_1} \varphi &=& \lambda_1(p) m_1(x) \varphi^{p-1} \ &{\rm in} \ \Omega, \\ 
	\varphi &=& 0 \ &{\rm in} \ \ \mathbb{R}^{N}\setminus\Omega,
\end{array}
\right.
\]
	for more details see \cite{Ianni}.
	
	We can suppose that $0 \leq \varphi \leq 1$ on $\Omega$, otherwise $\tilde{\varphi} = \displaystyle\frac{\varphi}{\Vert \varphi \Vert_\infty}$ satisfies the equation and $0 \leq \tilde{\varphi} \leq 1$ on $\Omega$. It follows that, for any constant $0 < c < \lambda_1(p)^{\frac{-1}{p-1-\alpha_1}}$, we have $\lambda_1(p) (c\varphi)^{p-1} \leq (c\varphi)^{\alpha_1}$, since $\alpha_1 < p-1$. Thus,
	\[
	(-\Delta_p)^{s_1} (\varphi) = c^{p-1} (-\Delta_p)^{s_1} \varphi = \lambda_1(p) m_1(x) (c\varphi)^{p-1} \leq m_1(x) (c\varphi)^{\alpha_1} \ \ \text{in} \ \ \Omega
	\]
	in other words $c\varphi$ is a lower solution of the problem \eqref{P2}. By comparison principle, choosing $c$ sufficiently small we can assume that $c\varphi < u_2$. Let us define the sequence $(w_n)_{n \in \mathbb{N}}$ by, $w_0 = u_2$ and $w_n$ the unique solution of the problem 
	\[
		\left\{
		\begin{array}{rrll}
	(-\Delta_p)^{s_1} w_n &=& m_1(x) w_{n-1}^{\alpha_1} \ &{\rm in} \ \Omega, \\  w_n &=& 0 \ &{\rm in} \  \mathbb{R}^{N}\setminus\Omega.
\end{array}
\right.
\]
	
	Note that, $w_n$ is positive  and $c\varphi \leq w_n \leq u_2$ for all $n \in \mathbb{N}$. Using regularity and compactness results we can prove that $w_n \to u$ in $C^0(\overline{\Omega})$ and $u$ is a solution of the problem \eqref{P2} with $c\varphi \leq u \leq u_2$. Thus, $u$ satisfies the equation 
	\[
		\left\{
		\begin{array}{rrll}
			(-\Delta_p)^{s_1} u &=& a(x) \vert u\vert^{\alpha_1}v_1^{\beta_1} & {\rm in} \ \ \Omega,\\
			u &=& 0  & {\rm in} \ \ \mathbb{R}^{N}\setminus\Omega
		\end{array}
		\right.
	\]
	and by uniqueness, we have $u = J_1(v_1)= u_1$. Consequently $u_1 \leq u_2$. Analogously, we conclude that $J_2$ is nondecreasing.
	\qed
	
	\begin{lema}\label{Lema4}
		\begin{enumerate}
			\item[(i)] A couple $(u,v) \in X_+\backslash \{0\}$ is a weak solution of the problem \eqref{1} for some $(\lambda, \mu) \in (\mathbb{R}_+^*)^2$ if, and only if, $u,v \in \text{int}(X_+)$ and $J_1(v) = \lambda^{\frac{-1}{p-1-\alpha_1}}u$ and $J_2(u) =  \mu^{\frac{-1}{q-1-\alpha_2}}v$;
			\item[(ii)] For any $\rho, \sigma > 0$ we have
			\begin{align*}
				(J_1 \circ J_2)(\rho u) &= \rho^{\frac{\beta_1\beta_2}{(p-1-\alpha_1)(q-1-\alpha_2)}} (J_1 \circ J_2)( u)\\
				(J_2 \circ J_1)(\sigma v) &= \sigma^{\frac{\beta_1\beta_2}{(p-1-\alpha_1)(q-1-\alpha_2)}} (J_2 \circ J_1)(v).
			\end{align*}
		\end{enumerate}
	\end{lema}
	\n {\bf Proof.} $(i)$ Indeed, if a couple $(u,v) \in X_+\backslash \{0\}$ is a weak solution of the problem \eqref{1} for some $(\lambda, \mu) \in (\mathbb{R}_+^*)^2$ then by strong maximum principle and Hopf's lemma, see \cite{Pezzo}, results that $u,v \in \text{int}(X_+)$. Moreover, since $v > 0$ in $\Omega$ and $(u,v)$ satisfies the problem \eqref{1}, we have, $u = 0 \  {\rm in} \ \ \mathbb{R}^{N}\setminus\Omega$ and 
	\begin{align*}
		(-\Delta_p)^{s_1} u &= \lambda a(x) \vert u \vert^{\alpha_1} v^{\beta_1 - 1 } v \ \ \text{in} \ \ \Omega. 
	\end{align*}
	Thus, for any $k>0$,
	\begin{align*}
		(-\Delta_p)^{s_1} (ku) &= k^{p-1}\lambda a(x) \vert u \vert^{\alpha_1} v^{\beta_1 } \ \ \text{in} \ \ \Omega \\
		&= k^{p-1-\alpha_1}\lambda a(x) \vert ku \vert^{\alpha_1} v^{\beta_1} \ \ \text{in} \ \ \Omega.
	\end{align*}
	
	Consequently $J_1(v) = ku$ if, and only if, 
	\[
 k= \lambda^{\frac{-1}{p-1-\alpha_1}}
	\]
	and in that case $J_1(v) = \lambda^{\frac{-1}{p-1-\alpha_1}}u$.
	
	Similarly, using the $u >0$, $v = 0 \  {\rm in} \ \ \mathbb{R}^{N}\setminus\Omega$ and 
	\begin{align*}
		(-\Delta_p)^{s_2} v &= \lambda a(x) \vert v \vert^{\alpha_2} u^{\beta_2 - 1 } u \ \ \text{in} \ \ \Omega 
	\end{align*}
	results that $J_2(u) = \mu^{\frac{-1}{q-1-\alpha_2}}v.$
	
	$(ii)$ Notice that, if we denote $J_2(\rho u) = \tilde{v}$ then in $\Omega$ we have,
	\[
	(-\Delta_q)^{s_2} \tilde{v} = a(x) \vert \tilde{v} \vert^{\alpha_2} (\rho u)^{\beta_2} = \rho^{\beta_2}  a(x) \vert \tilde{v} \vert^{\alpha_2} u^{\beta_2} 
	\]
	thus, for any $k>0$,
	\begin{align*}
		(-\Delta_q)^{s_2} (k\tilde{v}) &= k^{q-1} \rho^{-\beta_2} a(x) \vert \tilde{v} \vert^{\alpha_2} u^{\beta_2} \ \ \text{in} \ \ \Omega \\
		&= k^{q-1-\alpha_2}\rho^{-\beta_2} a(x) \vert k\tilde{v} \vert^{\alpha_2} u^{\beta_2} \ \ \text{in} \ \ \Omega.
	\end{align*}
	
	Consequently $J_2(\rho u) = \tilde{v} = kJ_2(u)$ if, and only if, 
	\[
 k= \rho^{\frac{\beta_2}{q-1-\alpha_2}}
	\]
	and in that case $J_2(\rho u) = \rho^{\frac{\beta_2}{q-1-\alpha_2}} J_2(u)$. Similarly, for any $\sigma >0$ results $J_1(\sigma v) = \sigma^{\frac{\beta_1}{p-1-\alpha_1}} J_1(v)$.
	
	Thus, for any $\rho, \sigma > 0$ we have
	\begin{align*}
		(J_1 \circ J_2)(\rho u) &= J_1(J_2(\rho u)) = J_1(\rho^{\frac{\beta_2}{q-1-\alpha_2}} J_2(u)) = \rho^{\frac{\beta_1\beta_2}{(p-1-\alpha_1)(q-1-\alpha_2)}} (J_1 \circ J_2)(u),\\
		(J_2 \circ J_1)(\sigma v) &= J_2(J_1(\sigma v)) = J_2(\sigma^{\frac{\beta_1}{p-1-\alpha_1}} J_1(u)) = \sigma^{\frac{\beta_1\beta_2}{(p-1-\alpha_1)(q-1-\alpha_2)}} (J_2 \circ J_1)(v).
	\end{align*}
	\qed
	
{\bf Proof of Theorem \ref{teor2}:}\\
\mbox{}

	It follows from Lemma \ref{Lema4} and condition \eqref{cond1} that
	\[
	(J_1 \circ J_2)(tu) = t(J_1 \circ J_2)(u) \ \ \text{and} \ \ (J_2 \circ J_1)(tv) = t(J_2 \circ J_1)(v),
	\]
	that is, the mappings $J_1 \circ J_2$ and $J_2 \circ J_1$ are homogeneous.
	
	Furthermore, by the strong comparison principle, both mappings are strongly monotone . The regularity results quoted before imply that $J_1 \circ J_2$ and $J_2 \circ J_1$ maps bounded sets into sets with compact closure.
	By Krein-Rutman Theorem, there exists a unique number $\Lambda_1 \in \mathbb{R}_+^{*}$ such that $(J_1 \circ J_2)(u_0) = \Lambda_1u_0$ holds for some $u_0 \in \left(C^0_{s_1}(\overline{\Omega})\right)_+$  which is unique up to a positive constant multiple. Analogously, there is a unique number $\Theta_1 \in \mathbb{R}_+^{*}$ such that $(J_2 \circ J_1) (v_0) = \Theta_1 v_0$ holds for some $v_0 \in\left(C^0_{s_2}(\overline{\Omega})\right)_+$.
	
	The $\displaystyle\frac{\beta_2}{q-1-\alpha_2}$-homogeneity of $J_2$ applied to $(J_1 \circ J_2)(u_0)$ results
	\[
	(J_2 \circ J_1)(J_2(u_0)) = J_2(\Lambda_1 u_0) = \Lambda_1^{\frac{\beta_2}{q-1-\alpha_2}} J_2(u_0).
	\]
	
	Analogously, we have
	\[
	(J_1 \circ J_2)(J_1(v_0)) = J_1(\Theta_1 v_0) = \Theta_1^{\frac{\beta_1}{p-1-\alpha_1}} J_1(v_0).
	\]
	
	The uniqueness of $\Theta_1$ and $v_0$ results 
	\[
	\Theta_1 = \Lambda_1^{\frac{\beta_2}{q-1-\alpha_2}}  \ \ \text{and} \ \ v_0 = k J_2(u_0) \ \ \text{for some} \ \ k> 0.
	\]
	
	On the other hand, the uniqueness of $\Lambda_1$ and $u_0$ provides us $\Lambda_1 = \Theta_1^{\frac{\beta_2}{q-1-\alpha_2}}$ and
	\[
	J_1(v_0) = J_1(kJ_2(u_0)) = k^{\frac{\beta_1}{p-1-\alpha_1}} (J_1 \circ J_2)(u_0) = k^{\frac{\beta_1}{p-1-\alpha_1}} \Lambda_1 u_0
	\]
	thus,
	\[
	\Lambda_1 = \Theta_1^{\frac{\beta_2}{q-1-\alpha_2}} \ \ \text{and} \ \ u_0 = k^{\frac{-\beta_1}{p-1-\alpha_1}} \Lambda_1^{-1} J_1(v_0).
	\]
	
	By Lemma \ref{Lema4}$(i)$, the pair $(u_0, J_2(u_0))$ solves \eqref{1} with $\lambda = \Lambda_1^{\alpha_1 - p + 1}$ and $\mu = 1$, since 
	\[
	J_2(u_0) = 1.J_2(u_0) \ \ \text{and} \ \ J_1(J_2(u_0)) = (J_1 \circ J_2)(u_0) = \Lambda_1 u_0 = \left(\Lambda_1^{\alpha_1 - p+1} \right)^{\frac{-1}{p-1-\alpha_1}} u_0.
	\]
	
	Analogously, $(J_1(v_0),v_0)$ solves \eqref{1} with $\lambda=1$ and $\mu = \Theta^{\alpha_2-p+1} = \Lambda_1^{-\beta_2}$, using the fact $\Lambda_1 = \Theta_1^{\frac{\beta_2}{q-1-\alpha_2}}$.
	
	Note that, if $\lambda$ and $\mu$ are positive numbers satisfying 
	\[
	\mu^{-\beta_2} \lambda^{\frac{-1}{p - 1 - \alpha_1}} = \Lambda_1
	\]
	then the pair $\rho = \lambda^{\frac{1}{p - 1 - \alpha_1}}$ and $\sigma_1 = 1$ satisfies the condition $(iii)$ of the Lemma \ref{Lema1} and thus $(\lambda, \mu)$ is a principal eigenvalue of \eqref{1}.
	
	Reciprocally, if $u,v \in X_+$ solves \eqref{1} for some $(\lambda,\mu)$ then, by Lemma \ref{Lema4} $(i)$
	\begin{align*}
		(J_1 \circ J_2)(u) &= J_1(J_2(u)) = J_1(\mu^{\frac{-1}{q-1-\alpha_2}}v) \\ 
		&= \left(\mu^{\frac{-1}{q-1-\alpha_2}} \right)^{\frac{\beta_1}{p-1-\alpha_1}}J_1(v) \\
		&= \mu^{\frac{-\beta_1}{(q-1-\alpha_2)(p-1-\alpha_1)}} \lambda^{\frac{-1}{p-1-\alpha_1}} u \\
		&= \mu^{\frac{-1}{\beta_2}} \lambda^{\frac{-1}{p-1-\alpha_1}} u.
	\end{align*}
	
	From uniqueness results of the Krein-Rutman Theorem, we obtain
	\[
	\mu^{-\beta_2} \lambda^{\frac{-1}{p - 1 - \alpha_1}} = \Lambda_1
	\]
	and $u= \rho u_0$ for some $\rho > 0$. Raising both sides of the above equation to the
	power $-\sqrt{\beta_2/(q-1-\alpha_2)}$, we have $(\lambda,\mu)$ satisfies \eqref{cond2} with $\Lambda_0 = \Lambda_1^{\sqrt{\frac{\beta_2}{q-1-\alpha_2}}}$. 
	Using the Lemma \ref{Lema1} $(iii)$, we have $v = \mu ^\frac{1}{\beta_2}\rho v_0$.
	\qed

\section{Proof of Theorem \ref{teor3}:}
\mbox{}

 Let us also denote by $(\varphi_\lambda, \psi_\mu)$ the positive eigenfunction associated to $(\lambda,\mu) \in \mathcal{C}_1$ with $\Vert \varphi_\lambda \Vert_{L^\infty(\Omega)} = 1$.

Properties $(P1)$ and $(P2)$ are immediate consequences of the Theorem \ref{teor2}. For property $(P3)$, assume that $(u,v) \notin X_+$, then by Hopf's lemma there is some $\gamma \in \mathbb{R}_+$ such that
	\[
	-u \geq \gamma\varphi_\lambda  \ \ \text{and} \ \ -v \geq \gamma^{\omega}\psi_\mu
	\]
	where $\omega = \frac{\beta_2}{q-1}$. Let $\overline{\gamma}$ be the minimum of such $\gamma's$. Suppose that $\overline{\gamma} > 0$ and $-u = \overline{\gamma} \varphi_\lambda$ in $\Omega$. Thus, since $\alpha_1 = \alpha_2 = 0$ we have of the second equation of \eqref{1}
	\[
	\mu b(x)  \vert -u \vert^{\beta_2 -1} (-u) = \mu b(x) \vert \overline{\gamma} \varphi_\lambda \vert^{\beta_2 - 1} \overline{\gamma}\varphi_\lambda
	\]
	and, consequently $-v = \overline{\gamma}^{\omega} \psi_{\mu}$ by uniqueness.
	
	Then, it remains to treat the case $-u \not\equiv \overline{\gamma} \varphi_{\lambda}$. Therefore, we have in $\Omega$,
	\begin{align*}
		(-\Delta_p)^{s_1} (\overline{\gamma}\varphi_{\lambda}) &= \overline{\gamma}^{p-1} (-\Delta_p)^s(\varphi_\lambda)\\
		&= \overline{\gamma}^{p-1}\lambda a(x) \vert \psi_\mu \vert^{\beta_1 - 1} \psi_\mu\\
		&= \lambda a(x) \vert \overline{\gamma}^{\frac{p-1}{\beta_1}}\psi_\mu \vert^{\beta_1} \\
		&= \lambda a(x) \vert \overline{\gamma}^{\frac{\beta_2}{q-1}}\psi_\mu \vert^{\beta_1} \\
		&= \lambda a(x) \vert \overline{\gamma}^{\omega}\psi_\mu \vert^{\beta_1-1} (\overline{\gamma}^{\omega} \psi_\mu)\\ 
		&\geq \lambda a(x) \vert -v \vert^{\beta_1-1} (-v) = (-\Delta_p)^{s_1}(-u).
	\end{align*}
	above we use $\psi_\mu > 0$ in $\Omega$ and the conditions $\alpha_1 = \alpha_2 =0$ and $\beta_1\beta_2 = (p-1)(q-1)$. Similarly we find 
	\[
	(-\Delta_q)^{s_2} (-v) \leq (-\Delta_q)^{s_2} (\overline{\gamma}^{\omega} \psi_\mu) \ \ \text{in} \ \ \Omega.
	\]
	
	Moreover, observe that $-u = \overline{\gamma} \varphi_\lambda = -v = \overline{\gamma}^{\omega} \psi_\mu = 0$ in $\mathbb{R}^N \backslash \Omega$. Since we assume that $-u \not\equiv \overline{\gamma}\varphi_\lambda$ and $-v \not\equiv \overline{\gamma}^{\omega}\psi_\mu$, by the
	strong comparison principle (Theorem 2.7, \cite{Mosconi}), $-u < \overline{\gamma} \varphi_{\lambda}$ and thus can find $0 < \varepsilon < 1$ such that, $-u \leq \varepsilon\overline{\gamma}\varphi_\lambda$ and $-v \leq (\varepsilon\overline{\gamma}^{\omega})\psi_\mu$, a contradiction with our definition of $\overline{\gamma}$. Therefore, $\overline{\gamma} =0$ and consequently we have the result. \qed
	
	\section{Proof of Theorem \ref{lower}} 
\mbox{}

Let $(\lambda_1, \mu_1)\in \mathcal{C}_1$ and $(\varphi,\psi) \in \text{int}(X_+)$ be a positive eigenfunction of the system \eqref{1} associated to $(\lambda_1, \mu_1)$. By applying the $L^{\infty}$-bound (\ref{abp}) to the first equation of \eqref{1}, we obtain
\begin{eqnarray}\label{i1}
||\varphi||_{L^\infty(\Omega)} &=& \sup_\Omega \varphi \leq C_{s_1,p}d^{\frac{ps_1}{p-1}}\lambda_1^{\frac{1}{p-1}} \|a\|^{\frac{1}{p-1}}_{L^\infty(\Omega)}||\varphi||^{\frac{\alpha_1}{p-1}}_{L^\infty(\Omega)} \|\psi\|^{\frac{\beta_1}{p-1}}_{L^\infty(\Omega)}
\end{eqnarray}
and to the second equation of \eqref{1}, we get
\begin{eqnarray}\label{i2}
||\psi||_{L^\infty(\Omega)} \leq C_{s_2,q}d^{\frac{qs_2}{q-1}}\mu_1^{\frac{1}{q-1}} \|b\|^{\frac{1}{q-1}}_{L^\infty(\Omega)} \|\varphi\|^{\frac{\beta_2}{q-1}}_{L^\infty(\Omega)}||\psi||^{\frac{\alpha_2}{q-1}}_{L^\infty(\Omega)} .
\end{eqnarray}

Thus, joining inequalities (\ref{i1}) and (\ref{i2}) and using the hypothesis 
\[
\beta_1 \beta_2 = (p-1-\alpha_1)(q-1-\alpha_2)\ \text{ and }\ \lambda_1^{\frac{1}{\theta}}\mu_1^{\frac{1}{\zeta}} = \Lambda_0,
\]
we derive \eqref{lb1}. This concludes the desired proof. \qed
	
\section{Proof of Theorem \ref{MP}}
\mbox{}

Note that it suffices to prove only that (i) $\Leftrightarrow$ (ii). In this case, we have (ii) $\Leftrightarrow$ (iii). In fact, we clearly have {\bf (SMP)} in $\Omega$ implies {\bf (WMP)} in $\Omega$. Conversely, suppose that {\bf (WMP)} holds in $\Omega$ and let $(u,v)$ be a weak solution of the system (\ref{supersol}). Then, $u, v \geq 0$ in $\Omega$ and, by (i), we get $\lambda, \mu \geq 0$. Thus, since $(-\Delta_m)^{s}$ satisfies the strong maximum principle (see \cite{Pezzo}), the conclusion of {\bf (SMP)} follows.

In order to show that {\bf (WMP)} in $\Omega$ leads to $(\lambda, \mu) \in \overline{\mathcal{R}_1} \setminus \mathcal{C}_1$, suppose instead that $(\lambda, \mu) \not\in \overline{\mathcal{R}_1} \setminus \mathcal{C}_1$. Let $(\lambda, \mu) \in \mathcal{C}_1$ and $(\tilde{\varphi},\tilde{\psi})$ be a positive eigenfunction corresponding to $(\lambda, \mu)$. Thus, $(-\tilde{\varphi},-\tilde{\psi})$ is a negative eigenfunction corresponding to $(\lambda, \mu)$. Therefore, {\bf (WMP)} fails in $\Omega$.

Suppose now that $(\lambda, \mu) \in \mathbb{R}^2$ is a fixed couple outside of $\overline{\mathcal{R}_1}$. If $(\lambda, \mu) \in (\mathbb{R}_+^*)^2$, we obtain $\lambda > \lambda_1$ and $\mu > \mu_1$, where $(\lambda_1, \mu_1)$ is a principal eigenvalue of (\ref{1.3}) with $\frac{\mu}{\lambda}=\frac{\mu_1}{\lambda_1}$. We denote by $(\varphi,\psi)$ a positive eigenfunction associated to $(\lambda_1, \mu_1)$. Then, $(-\varphi,-\psi)$ satisfies

\[
\left\{
\begin{array}{llll}
(-\Delta_p)^{s_1} (-\varphi) - \lambda a(x) \vert -\psi\vert^{\beta_1-1}(-\psi) &=& - \lambda_1 a(x) \psi^{\beta_1} + \lambda a(x) \psi^{\beta_1}\\
&=& (\lambda - \lambda_1)a(x) \psi^{\beta_1} \geq (\not\equiv)\ 0  \ \ & {\rm in}\ \Omega,\\
(-\Delta_q)^{s_2} (-\psi) - \mu b(x)\vert -\varphi\vert^{\beta_2-1}(-\varphi) &=& -\mu_1 b(x) \varphi^{\beta_2} + \mu b(x) \varphi^{\beta_2} \\
&=& (\mu - \mu_1)b(x) \varphi^{\beta_2} \geq (\not\equiv)\ 0 \ \ & {\rm in}\ \Omega
\end{array}
\right. 
\]
in sense weak and $-\varphi = 0 =-\psi$ in $\mathbb{R}^N\setminus\Omega$. Since, $-\varphi,-\psi<0$ in $\Omega$, we have {\bf (WMP)} fails in $\Omega$.

Now, assume that $\lambda < 0$. Thus, there is $(\lambda_1, \mu_1)\in \mathcal{C}_1$ with $\lambda_1 > 0$ small enough (and so $\mu_1>0$ large enough) so that $\lambda < - \lambda_1$ and $\mu > - \mu_1$. Then, $(-\varphi,\psi)$ satisfies

\[
\left\{
\begin{array}{llll}
(-\Delta_p)^{s_1} (-\varphi) - \lambda a(x) \psi^{\beta_1} &=& -\lambda_1 a(x) \psi^{\beta_1} - \lambda a(x) \psi^{\beta_1} \\
&=& -(\lambda + \lambda_1)a(x) \psi^{\beta_1} \geq (\not\equiv)\ 0 \ \ & {\rm in}\ \Omega, \\
(-\Delta_q)^{s_2} \psi - \mu b(x)\vert -\varphi\vert^{\beta_2-1}(-\varphi) &=& \mu_1 b(x) \varphi^{\beta_2} + \mu b(x) \varphi^{\beta_2} \\
&=& (\mu+\mu_1) b(x) \varphi^{\beta_2} \geq (\not\equiv)\ 0\ \ & {\rm in}\ \Omega
\end{array}
\right. 
\]
in sense weak and $-\varphi = 0 =\psi$ in $\mathbb{R}^N\setminus\Omega$. But, $-\varphi<0$ in $\Omega$ and so {\bf (WMP)} fails in $\Omega$.

For the remaining case $\lambda \geq 0$ and $\mu < 0$, there is $(\lambda_1, \mu_1)\in \mathcal{C}_1$ with $\lambda_1 > 0$ large enough (and so $\mu_1>0$ small enough) so that $\lambda > -\lambda_1$ and $\mu < -\mu_1$. Therefore, $(\varphi,-\psi)$ satisfies

\[
\left\{
\begin{array}{llll}
(-\Delta_p)^{s_1} \varphi - \lambda a(x) \vert -\psi\vert^{\beta_1-1}(-\psi) &=& \lambda_1 a(x) \psi^{\beta_1} + \lambda a(x)\psi^{\beta_1} \\
&=& (\lambda + \lambda_1)a(x) \psi^{\beta_1} \geq (\not\equiv)\ 0\ \ & {\rm in}\ \Omega, \\
(-\Delta_q)^{s_2} (-\psi) - \mu b(x) \varphi^{\beta_2} &=& - \mu_1 b(x) \varphi^{\beta_2} - \mu b(x) \varphi^{\beta_2} \\
&=& -(\mu + \mu_1)b(x) \varphi^{\beta_2} \geq (\not\equiv)\ 0\ \ & {\rm in}\ \Omega
\end{array}
\right. 
\]
in sense weak and $\varphi = 0 =-\psi$ in $\mathbb{R}^N\setminus\Omega$. Therefore, $-\psi<0$ in $\Omega$ and so again {\bf (WMP)} fails in $\Omega$.

Conversely, we next prove that {\bf (WMP)} holds in $\Omega$ for any couple $(\lambda, \mu) \in \overline{\mathcal{R}_1} \setminus \mathcal{C}_1$. Since $(-\Delta_m)^{s}$ satisfies weak maximum principle in $\Omega$, we obtain {\bf (WMP)} holds in $\Omega$ if either $\lambda = 0$ and $\mu \geq 0$ or $\lambda \geq 0$ and $\mu = 0$. Let $(\lambda,\mu) \in \mathcal{R}_1$ and $(u,v)$ be a weak solution of the system (\ref{supersol}).  Notice that $\lambda < \lambda_1$ and $\mu < \mu_1$, where $(\lambda_1, \mu_1)$ is a principal eigenvalue of (\ref{1.3}) with $\frac{\mu}{\lambda}=\frac{\mu_1}{\lambda_1}$. Let $(\varphi,\psi)$ be a positive eigenfunction corresponding to $(\lambda_1, \mu_1)$.

We affirm that $u,v\geq 0$ in $\Omega$. In fact, suppose by contradiction that $u$ or $v$ is negative somewhere in $\Omega$. Thus, by Hopf's lemma and strong comparison principle for the fractional $m$-Laplacian (see \cite{CFGQ, Pezzo}), there is some $\gamma>0$ such that 
\[
-u\leq\gamma\varphi\ \text{and}\ -v\leq\gamma^\omega\psi\ \ \text{in}\ \ \Omega,
\]
where $\omega=\frac{p-1}{\beta_1}$. Let $\overline{\gamma}$ be the minimum of such $\gamma's$. Then, $\overline{\gamma}>0$. Since $\lambda < \lambda_1$ and $\mu < \mu_1$, we obtain
\[
\left\{
\begin{array}{llll}
(-\Delta_p)^{s_1} (-u) \leq \lambda a(x) \vert -v\vert^{\beta_1-1}(-v)\leq  \lambda a(x) (\overline{\gamma}^\omega\psi)^{\beta_1}\leq (\not\equiv)\ \lambda_1 a(x)  (\overline{\gamma}^\omega\psi)^{\beta_1}=(-\Delta_p)^{s_1}(\overline{\gamma}\varphi), \\
(-\Delta_q)^{s_2} (-v) \leq \mu b(x) \vert -u\vert^{\beta_2-1}(-u)\leq  \mu b(x) (\overline{\gamma}\varphi)^{\beta_2}\leq (\not\equiv)\ \mu_1 b(x)  (\overline{\gamma}\varphi)^{\beta_2}=(-\Delta_q)^{s_2}(\overline{\gamma}^\omega\psi)  
\end{array}
\right.
\]
in $\Omega$ in sense weak with $\overline{\gamma}\varphi\geq -u$ and $\overline{\gamma}^\omega\psi\geq - v$ in $\mathbb{R}^N\setminus\Omega$. So, using (\ref{h1}) and (\ref{h2}), it follows from the weak and strong comparison principles to each above equation (see \cite{CFGQ}) that $ -u< \overline{\gamma}\varphi$ and $-v<\overline{\gamma}^\omega\psi$ in $\Omega$. Then, we can find $0<\varepsilon <1$ such that $-u\leq\varepsilon\overline{\gamma}\varphi$ and $-v\leq(\varepsilon\overline{\gamma})^\omega\psi$ in $\Omega$, contradicting the definition of $\overline{\gamma}$. Therefore, $u, v \geq 0$ in $\Omega$. This concludes the desired proof. \qed

\section{Proof of Theorem \ref{CP}}
\mbox{}

Using (\ref{h3}), (\ref{h4}) and weak and strong comparison principles (see \cite{CFGQ}) and arguing in an analogous way as in the proof of Theorem \ref{MP}, we see that it suffices to prove only that (i) $\Leftrightarrow$ (ii).

(ii) $\Rightarrow$ (i) Taking $u, v \equiv 0$ in $\Omega$ and using that $(z,w)$ satisfies (\ref{h1}) and (\ref{h2}), we get {\bf (WCP)} in $\Omega$ implies {\bf (WMP)} in $\Omega$. So, by Theorem \ref{MP}, {\bf (WCP)} in $\Omega$ implies $(\lambda, \mu) \in \overline{\mathcal{R}_1} \setminus \mathcal{C}_1$. 

(i) $\Rightarrow$ (ii) We consider $(\lambda, \mu) \in \overline{\mathcal{R}_1} \setminus \mathcal{C}_1$. By weak maximum and comparison principles in $\Omega$, the conclusion is direct in the cases that $\lambda = 0$ or $\mu = 0$. Then, it suffices to consider $(\lambda, \mu) \in \mathcal{R}_1$. Set
\begin{eqnarray*}
&& f_1(x):= (-\Delta_p)^{s_1} u - \lambda a(x) v^{\beta_1}, \ \ f_2(x):= (-\Delta_p)^{s_1} z - \lambda a(x) w^{\beta_1},\\
&& g_1(x):= (-\Delta_q)^{s_2} v - \mu b(x) u^{\beta_2}, \ \ g_2(x):= (-\Delta_q)^{s_2} w - \mu b(x) z^{\beta_2}.
\end{eqnarray*}

Note that $0\leq f_1\leq f_2$ and $0\leq g_1\leq g_2$ in $\Omega$. Since $(u,v)$ satisfies (\ref{h1}) and (\ref{h2}) and is a weak solution of the problem \eqref{supersol}, by Theorem \ref{MP}, if $f_1,g_1\equiv 0$ in $\Omega$ we obtain $u, v \equiv 0$ in $\Omega$ and if $f_1+g_1\not\equiv 0$ in $\Omega$, we have $u, v > 0$ in $\Omega$. Notice that, when $f_1, g_1 \equiv 0$ in $\Omega$, the conclusion follows readily from {\bf (WMP)} in $\Omega$. Suppose then $f_1+g_1\not\equiv 0$ in $\Omega$ (and so $u, v > 0$ in $\Omega$). Then, since $(z,w)$ is also a weak solution of \eqref{supersol}, by {\bf (SMP)}, we derive $z, w > 0$ in $\Omega$. 

To conclude the proof of {\bf (WCP)}, it suffices to prove that $u \leq z$ in $\Omega$. This follows directly by using the fact that $(-\Delta_q)^{s_2}$ satisfies weak comparison principle in $\Omega$ (see \cite{CFGQ}). Arguing by contradiction, suppose that $u > z$ somewhere in $\Omega$. In this case, the set $\Gamma := \{\gamma>0 ; z > \gamma u\ {\rm and} \ w > \gamma^{\omega} v\ {\rm in} \ \Omega\}$, where $\omega=\frac{p-1}{\beta_1}$, is nonempty by Hopf's Lemma (see \cite{Pezzo}) and is also upper bounded. Set $\overline{\gamma}:=\sup \Gamma > 0$. Note that $z \geq \overline{\gamma} u$ and $w \geq \overline{\gamma}^{\omega} v$ in $\Omega$. Note also that $\overline{\gamma} < 1$. Now, since $f_1+g_1\not\equiv 0$ and $f_2+g_2\not\equiv 0$ in $\Omega$, we derive

\[
\left\{
\begin{array}{lll}
(-\Delta_p)^{s_1} (\overline{\gamma}u) =  \lambda a(x) (\overline{\gamma}^\omega v)^{\beta_1}+\overline{\gamma}^{p-1}f_1\leq (\not\equiv)\ \lambda a(x)  w^{\beta_1}+f_2=(-\Delta_p)^{s_1}(z) \ \ & {\rm in}\ \Omega, \\
(-\Delta_q)^{s_2} (\overline{\gamma}^\omega v) =  \mu b(x) (\overline{\gamma}u)^{\beta_2}+\overline{\gamma}^{\beta_2}g_1\leq (\not\equiv)\ \mu b(x)  z^{\beta_2}+g_2=(-\Delta_q)^{s_2}(w) \ \ & {\rm in}\ \Omega
\end{array}
\right. 
\]
in sense weak with $\overline{\gamma}u\leq z$ and $\overline{\gamma}^\omega v\leq w$ in $\mathbb{R}^N\setminus\Omega$. Then, using (\ref{h3}) and (\ref{h4}), it follows from the weak and strong comparison principles to each above equation (see \cite{CFGQ}) that $ z> \overline{\gamma} u$ and $w>\overline{\gamma}^\omega v$ in $\Omega$. Therefore, we can find $0<\varepsilon <1$ such that $z \geq (\overline{\gamma} + \varepsilon) u$ and $w \geq (\overline{\gamma}+ \varepsilon)^\omega v$ in $\Omega$. But this contradicts the definition of $\overline{\gamma}$. Hence, we complete the wished proof of theorem. \qed

\section{Proof of Theorem \ref{sm}}
\mbox{}

By Theorem \ref{MP}, the necessity of hypotheses $\lambda,\mu\geq 0$ and the equivalence between (ii) and (iii) follow immediately. Then, it suffices to show that the item (i) implies (ii).

We suppose that $\lambda \geq 0$ and $\mu \geq 0$. If either $\lambda = 0$ or $\mu = 0$, then by Theorem \ref{MP}, the desired {\bf (WMP)} follow.

Finally, suppose that $\lambda > 0$ and $\mu > 0$. Consider now the constant $\eta>0$ defined by

\[
\eta = \frac{1}{\left[\lambda^{\frac{1}{\theta}}\mu^{\frac{1}{\zeta}}C_{s_1,p}^{\frac{p-1}{\theta}}C_{s_2,q}^{\frac{q-1}{\zeta}} \|a\|^{\frac{1}{\theta}}_{L^\infty(\Omega)} \|b\|^{\frac{1}{\zeta}}_{L^\infty(\Omega)}\right]^{\frac{\beta_1\beta_2}{ps_1\theta+qs_2\zeta}}},
\]
where $d=\operatorname{diam }(\Omega)$, $\theta=\sqrt{\beta_1(p-1)}$, $\zeta=\sqrt{\beta_2(q-1)}$ and $C_{s_1,p}$ and $C_{s_2,q}$ are the explicit constants of $L^{\infty}$-bound (\ref{abp}). Thus, by using the lower estimate \eqref{lb1}, we get

\[
\Lambda_0 \geq \frac{1}{C_{s_1,p}^{\frac{p-1}{\theta}}C_{s_2,q}^{\frac{q-1}{\zeta}}d^{\frac{ps_1}{\theta}+\frac{qs_2}{\zeta}} \|a\|^{\frac{1}{\theta}}_{L^\infty(\Omega)} \|b\|^{\frac{1}{\zeta}}_{L^\infty(\Omega)}}>\lambda^{\frac{1}{\theta}}\mu^{\frac{1}{\zeta}}
\]
whenever $d < \eta$. Then, we derive $(\lambda, \mu) \in \mathcal{R}_1$ for such domains and so, by Theorem \ref{MP} the assertion (ii) holds. This concludes the proof. \qed

	\subsection*{Acknowledgement}
	
	The first author was partially supported by FAPEMIG/ APQ-02375-21, APQ-04528-22, FAPEMIG/RED-00133-21 and CNPq Process 307575/2019-5.\\
	The second author was partially supported by CNPq/Brazil (PQ 316526/2021-5) and Fapemig/Brazil (Universal-APQ-00709-18).

	\bibliographystyle{amsplain}

\end{document}